\documentclass[11pt]{amsart}
\usepackage{amssymb,amsmath,amsthm}
\usepackage[hmargin=3cm,vmargin=3.5cm]{geometry}
\usepackage{graphicx}
\usepackage{amscd}
\usepackage[all, knot]{xy}
\usepackage{epsfig}
\usepackage{psfrag}
\usepackage{lscape}
\newtheorem{theorem}{Theorem}[section]

\newtheorem{definition}[theorem]{Definition}

\title{Almost alternating knots with 12 crossings and Turaev genus}

\author{SLAVIK JABLAN}

\begin{document}
\maketitle

\begin{abstract}
Among non-alternating knots with $n\le 12$ crossings given in
\cite{1} Turaev genus is not known for 191 knot. For 154 of them we
show that they are almost alternating, so their Turaev genus is 1.
\end{abstract}

Almost alternating knots are defined in \cite{2}:

\begin{definition}
A diagram $D$ of a knot $K$ in the 3-sphere is almost alternating if
one crossing change makes the diagram alternating. A knot $K$ is
almost alternating if $K$ has an almost alternating diagram and $K$
does not have an alternating diagram.
\end{definition}

In \cite{2} the authors proved that all knots with at most $n=11$
crossings, except knots $K11n95=2\,0.-2\,-1.-2\,0.2$, $K11n118=2\,
0.-3.-2\,0.2$, and $K11n183=9^*.-2:.-2$ are almost alternating. For
the knot $K11n183$ with the minimal DT code
$\{\{11\},\{6,-14,10,-18,2,$ $22,20,-4,12,-8,16\}\}$, in the book
\cite{3} is shown that it is almost alternating, because it has
almost alternating representation $8^*2\,0.2\,0.-1.2\,0.2\,0$, with
the DT code $\{\{12\},\{-6,10,22,18,2,$ $16,24,20,8,12,4,14\}\}$,
and the status of the remaining two knots, $K11n95$ and $K11n118$ is
still undetermined.

To a knot diagram $D$, Turaev associated a closed orientable surface
embedded in $S^3$, called the Turaev surface. The Turaev genus
$g_T(K)$ of a knot $K$ is first defined in \cite{4}, as the minimal
number of the genera of the Turaev surfaces of diagrams of $K$. We
cite the description of Turaev surface and the definition of Turaev
genus from \cite{5}.

Every knot diagram $D$ has an associated Turaev surface $\Sigma D$
constructed as follows. Think of the diagram $D$ as a subset of
$S^2$, which itself is a subset of $S^3$. Replace the crossings of
$D$ with saddles so that the $A$-smoothing lies on one side of $S^2$
and the $B$-smoothing lies on the other side of $S^2$. Replace the
arcs of $D$ not near crossings with bands orthogonal to $S^2$. The
resulting surface is a cobordism between the all-$A$ Kauffman state
and the all-$B$ Kauffman state. The Turaev surface $\Sigma D$ is
obtained by capping off the boundary components of this cobordism
with disks.

\begin{figure}[th]
\centerline{\psfig{file=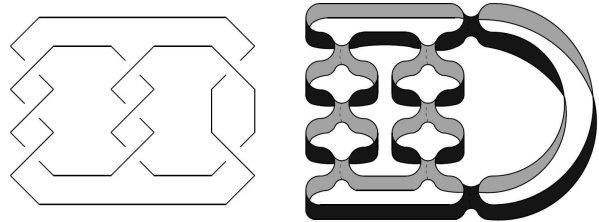,width=3.20in}} \vspace*{8pt}
\caption{Diagram of the knot $8_{19}=3,3,-2$ and its corresponding
Turaev surface \cite{5}. \label{ff}}
\end{figure}

\begin{definition}
The Turaev genus of a knot $K$, denoted $g_T(K)$, is defined as $g_T
(K) = min\{g(\Sigma D)\,\, |\,\, \text{ D is a diagram of K}\}$.
\end{definition}

For every alternating knot, its Turaev genus is 0, and for every
almost alternating knot it is 1 \cite{6}.

Since every almost alternating knot $K$ has the Turaev genus
$g_T(K)=1$, Turaev genus is known and equal to 1 for all
non-alternating knots with at most $n=11$ crossings, except $K11n95$
and $K11n118$, for which is unknown. These results are given in the
{\it Knot Info} {\tt (http://www.indiana.edu/$\sim $knotinfo/}) by
C. Livingston and J. C Cha \cite{1}.

All knots in this paper are given by their {\it Knotscape} symbol
and in the Conway notation, where for knots derived from basic
polyhedra is used notation from the program {\it LinKnot} {\tt
http://math.ict.edu.rs/}, which is used for all computations. Knots
are also given by their DT-codes.

For the computations we used two methods: replacement of rational
tangles with equivalent rational tangles with exactly one $-1$
elementary tangle, and the derivation of "source links" \footnote {A
link in the Conway notation containing only tangles $\pm 1$ and $\pm
2$ is called source link. Other links belonging to its family are
obtained by substituting twists $\pm 2$ by larger twists.} with
exactly one $-1$ tangle and their extension to larger almost
alternating knots.

Looking for Turaev genus of non-alternating knots with $n\le 12$
crossings in the {\it Knot Info} \cite{1} we conclude that it is
unknown for 191 knot. From them, we succeeded to show that 154 knots
are almost alternating by finding their explicit almost alternating
representations. These knots are given in the first two tables, and
for all of them Turaev genus is equal to 1. For the remaining 35
non-alternating knots with $n=12$ crossings their Turaev genus is
still unknown.

\bigskip
\noindent
\begin{tabular}{|c|c|c|} \hline
$K12n176$ & $.(2\,1,2).-2.2\,0$ & $.(2\,1,2).4\,-1\,1\,1.2\,0$ \\
\hline $K12n204$ & $.2.-2.(2,2\,1)\,0$ &
$.2.4\,-1\,1\,1.(2,2\,1)\,0$ \\ \hline $K12n257$ &
$2:(-2,-2\,-1)\,0:-2\,0$ & $(2,2\,1)\,0.-1.2.2:2\,0$ \\ \hline
$K12n258$ & $.-2\,-1\,-1.2\,0.2.2\,0$ & $.2\,2\,-1\,0.2\,0.2.2\,0$
\\ \hline $K12n281$ & $2\,1\,1\,0.2.-2.2\,0$ &
$2\,1\,1\,0.2.4\,-1\,1\,1.2\,0$ \\ \hline $K12n282$ &
$2.2\,0.2.-2\,-1\,-1\,0$ & $2.2\,0.2.2\,2\,-1$ \\ \hline $K12n286$ &
$.2\,1\,1\,1.-2\,0.2$ & $.2\,1\,1\,1.4\,-1\,1\,1\,0.2$ \\ \hline
$K12n299$ & $-2.-2\,0.-2.2\,1\,1\,0$ &
$8^*2.2.1.2\,0.2.2\,2\,0.1.-1$
\\ \hline $K12n316$ & $2\,0.2\,1.-2\,0.3\,0$ &
$2\,0.2\,1.4\,-1\,1\,1\,0.3\,0$ \\ \hline $K12n317$ &
$(2\,1\,1,-2\,-1)\,(3,2)$ & $(2\,1\,1,2\,2\,-1\,1\,0)\,(3,2)$ \\
\hline $K12n319$ & $2\,0.-2\,-1.-2\,0.-3\,0$ &
$8^*2\,1.2.-1.2\,0.1.1.3\,0.1$ \\ \hline $K12n321$ &
$3.2\,1\,0.-2.2\,0$ & $3.2\,1\,0.4\,-1\,1\,1.2\,0$ \\ \hline
$K12n322$ & $2\,3,2\,1\,1,-2\,-1$ & $2\,3,2\,1\,1,2\,2\,-1\,1\,0$ \\
\hline $K12n324$ & $-3.-2\,-1\,0.-2.2\,0$ &
$8^*2.2.2.2\,1\,0.2.2\,0.1.-1$ \\ \hline $K12n327$ &
$2\,0.2\,1.-2\,0.2\,1\,0$ & $2\,0.2\,1.4\,-1\,1\,1\,0.2\,1\,0$ \\
\hline $K12n349$ & $2\,1:-3\,0:-2\,-1\,0$ & $2\,1\,0.-1.3.3:2\,0$ \\
\hline $K12n353$ & $-2\,-1\,-1\,0:-2\,-1\,0:-2\,0$ &
$2\,1\,1:2\,-1:2\,1$ \\ \hline $K12n357$ & $8^*-2\,-1.-2\,0.-2$ &
$9^*2.2.2.2\,1.1.1.1.1.-1$ \\ \hline $K12n359$ & $3:2\,1:-2\,-1\,0$
& $3:2\,1:2\,2\,-1\,1$ \\ \hline $K12n360$ &
$-2\,-1.2\,0.-2.2.-2\,0$ & $102^*2.1.1.2\,0.1.-1.1.1.1.2$ \\ \hline
$K12n361$ & $8^*2\,0.2\,1:.-2$ & $8^*2\,0.2\,1:.4\,-1\,1\,1$ \\
\hline
$K12n376$ & $2.2.-2.-2\,0.-2\,-1$ & $9^*2.2.2\,0.1.2\,0.-1.1.1.2$ \\
\hline $K12n381$ & $-2\,-1.-2.-2.2\,1\,0$ &
$8^*2.2\,0.1.2\,1.2\,0.1.2\,0.-1$ \\ \hline $K12n382$ &
$-2\,-1.2.-2.2\,1\,0$ & $6^*2.2\,1.2\,1.2\,0.2\,0.-1$ \\ \hline
$K12n394$ & $8^*2.-2\,0:-2\,-1\,0$ &
$124^*2\,0.-1.1.1.1.2\,1.2.2.1.1.1.1$ \\ \hline $K12n395$ &
$-3\,0:2\,1\,0:2\,1\,0$ & $2\,1\,-1:2\,1\,0:2\,1\,0$ \\ \hline
$K12n398$ & $2\,1\,1:-2\,-1\,0:-2\,0$ & $2\,1\,1\,0.-1.2.2\,1:2\,0$
\\ \hline $K12n399$ & $-2\,-1\,-1:-2\,-1\,0:-2\,0$ &
$2\,1\,1\,0:2\,-1:2\,1$ \\ \hline $K12n403$ & $2.-2\,0.-2.-2\,-2\,0$
& $111^*1.1.2\,0.2.-1.1.1.2\,2\,0.1.1.2\,0$ \\ \hline $K12n417$ &
$2\,1:-2\,-1\,0:-2\,-1\,0$ & $2\,1\,0.-1.3.2\,1:2\,0$ \\ \hline
$K12n425$ & $8^*-2\,-1::2\,1$ & $8^*3\,-1\,0::2\,1$ \\ \hline
$K12n428$ & $-2\,-1.2.2.2\,0.2$ & $2\,2\,-1\,1\,0.2.2.2\,0.2$ \\
\hline $K12n430$ & $8^*2:-2\,-1\,0:-2\,0$ &
$15120^*1.1.1.-1.1.1.1.1.1.1.2\,0.1.1.1.1$ \\ \hline $K12n434$ &
$8^*-2\,-2.2\,0$ & $8^*2\,1\,1\,-1\,0.2\,0$ \\ \hline $K12n436$ &
$2.-2\,-1.-2\,0.2\,1\,0$ & $8^*2\,1\,0.2\,0.-1.2\,1$ \\ \hline
$K12n440$ & $8^*2\,1.-2:.2$ & $8^*2\,1.4\,-1\,1\,1:.2$ \\ \hline
$K12n444$ & $8^*2\,1.-1.-1.-1.3\,0$ & $9^*3\,0.-1.1.1.2.2\,0.1.1.1$ \\
\hline $K12n450$ & $8^*2.2\,1.-2\,0$ & $8^*2.2\,1.4\,-1\,1\,1\,0$ \\
\hline $K12n453$ & $8^*2:.2\,1\,0:.-2\,0$ &
$8^*2:.2\,1\,0:.4\,-1\,1\,1\,0$ \\ \hline $K12n459$ &
$8^*2\,1\,1\,0::-2\,0$ & $8^*2\,1\,1\,0::4\,-1\,1\,1\,0$ \\ \hline
$K12n463$ & $8^*-2\,0.2:2\,1\,0$ & $8^*4\,-1\,1\,1\,0.2:2\,1\,0$ \\
\hline $K12n465$ & $102^*2\,0:.-2$ & $102^*2\,0:.4\,-1\,1\,1$ \\
\hline $K12n481$ & $-2.-2\,0.-2.2\,2\,0$ &
$8^*2.2.1.2\,0.2.2\,1\,1\,0.1.-1$ \\ \hline $K12n485$ &
$-2\,-1.-2.-2.2.2$ & $101^*2.-1.2.1.2\,0.2.1.1.1.1$ \\ \hline
$K12n486$ & $2\,1\,0.-2.-2\,0.-3\,0$ & $8^*2\,1.2.-1.2.1.1.3.1$ \\
\hline $K12n489$ & $2\,1\,0.2.2\,0.-3\,0$ &
$2\,1\,0.2.2\,0.2\,1\,-1$ \\ \hline $K12n497$ & $2\,1\,2\,1,-3,2\,1$
& $2\,1\,2\,1,2\,1\,-1\,0,2\,1$ \\ \hline

\end{tabular}

\begin{tabular}{|c|c|c|} \hline

$K12n508$ & $(-2\,-1,2\,1)\,(2\,1,2\,1)$ &
$(2\,2\,-1\,1\,0,2\,1)\,(2\,1,2\,1)$
\\ \hline $K12n509$ & $-3:-2\,-1\,0:-2\,-1\,0$ &
$8^*2.-1.2.1.2\,0.3.2\,0.1$ \\ \hline $K12n510$ &
$-3:2\,1\,0:2\,1\,0$ & $2\,1\,-1\,0:2\,1\,0:2\,1\,0$ \\ \hline
$K12n519$ & $3:-2\,-1\,0:-2\,-1\,0$ & $3\,0.-1.3.2\,1:2\,0$ \\
\hline $K12n520$ & $2.2\,1.-3.2\,0$ & $2.2\,1.2\,1\,-1\,0.2\,0$ \\
\hline $K12n521$ & $8^*2.-2\,0.2\,1\,0$ &
$8^*2.4\,-1\,1\,1\,0.2\,1\,0$ \\ \hline $K12n524$ &
$-3\,0:-2\,-1\,0:-2\,-1\,0$ & $8^*2.-1.2.1.2\,0.3\,0.2\,0.1$ \\
\hline $K12n527$ & $2.2\,1.-2.2.2\,0$ & $2.2\,1.4\,-1\,1\,1.2.2\,0$
\\ \hline $K12n535$ & $8^*-2.-2\,0.2\,1\,0$ &
$8^*2.-1.3\,0.2\,0.1.2.1.1$ \\ \hline $K12n536$ &
$8^*-2.-2\,0.-2\,-1\,0$ & $9^*2\,1\,0.2\,0.1.1.1.1.1.-1.2\,0$ \\
\hline $K12n542$ & $-2.-2\,-1.-2.2.2\,0$ &
$9^*2\,1\,0.2\,0.1.1.2\,0.1.1.-1.1$ \\ \hline $K12n544$ &
$2.-2.-2\,-1.-2.2$ & $9^*2.2\,1.1.1.1.2.1.-1.1$ \\ \hline $K12n551$
& $-2\,-1\,-1\,-1\,0:-2\,0:-2\,0$ & $2\,1\,1\,1:2\,-1:2$ \\ \hline
$K12n552$ & $2.-2.-2\,0.2\,1\,1\,0$ & $8^*2\,1\,1\,0.2\,0.-1.2:1$ \\
\hline $K12n557$ & $-2\,-2\,0:-2\,-1\,0:-2\,0$ & $2\,2:2\,-1:2\,1$
\\ \hline $K12n559$ & $8^*3\,0:.-2:.2\,0$ &
$8^*3\,0:.4\,-1\,1\,1:.2\,0$ \\ \hline $K12n561$ & $8^*3:.-2:.2\,0$
& $8^*3:.4\,-1\,1\,1:.2\,0$ \\ \hline $K12n562$ & $.2\,1.-2.4\,0$ &
$.2\,1.4\,-1\,1\,1.4\,0$ \\ \hline $K12n563$ & $3:2\,2\,0:-2\,0$ &
$3:2\,2\,0:4\,-1\,1\,1\,0$ \\ \hline $K12n567$ & $9^*-2\,0.2\,1\,0$
& $9^*4\,-1\,1\,1\,0.2\,1\,0$ \\ \hline $K12n568$ &
$9^*2\,0.2::-2\,0$ & $9^*2\,0.2::4\,-1\,1\,1\,0$ \\ \hline $K12n569$
& $-2\,-2\,0:-3\,0:-2\,0$ & $8^*2.-1.2.1.2.2\,2\,0.1.1$ \\ \hline
$K12n579$ & $-2.2.-2.2\,1.2\,0$ & $8^*2.-1.2\,0.2\,0.1.2\,1\,0.1.1$
\\ \hline $K12n580$ & $-3:-2\,-2\,0:-2\,0$ & $8^*2.-1.2.1.1.3.2\,1.1$
\\ \hline $K12n584$ & $2\,1.-2.3.2$ & $2\,1.4\,-1\,1\,1.3.2$ \\
\hline $K12n589$ & $9^*2\,1\,0:-2\,0$ & $9^*2\,1\,0:4\,-1\,1\,1\,0$ \\
\hline $K12n594$ & $3.-2\,-1.-2\,0.2\,0$ & $8^*2\,0.2\,0.-1.2\,1:2$
\\ \hline $K12n612$ & $-2\,-2\,-1\,0:-2\,0:-2\,0$ &
$2\,2\,1:2\,-1:2$ \\ \hline $K12n617$ & $9^*-2\,-1\,0.-2$ &
$16229^*1.1.1.1.-1.1.1.1.1.1.1.1.1.1.1.1$ \\ \hline $K12n618$ &
$8^*2\,1\,0.-2\,0:2\,0$ & $8^*2\,1\,0.4\,-1\,1\,1\,0:2\,0$ \\ \hline
$K12n619$ & $8^*2\,1\,0:2\,0:-2\,0$ &
$8^*2\,1\,0:2\,0:4\,-1\,1\,1\,0$
\\ \hline $K12n624$ & $8^*-2\,-1\,0:.2:.-2\,0$ &
$8^*2.2\,0.2.-1.2\,0.2\,0.1.1$ \\ \hline $K12n629$ &
$8^*-2\,-1\,0:-2\,0:2\,0$ & $111^*1.2\,0.1.1.1.-1.1.3\,0.1.1.2$ \\
\hline $K12n630$ & $9^*.-2\,-1:.-2$ & $8^*3\,0.2\,0.-1.2\,0.2\,0$ \\
\hline $K12n631$ & $8^*2\,1\,0:-2\,0:2\,0$ &
$8^*2\,1\,0:4\,-1\,1\,1\,0:2\,0$ \\ \hline $K12n638$ &
$2:-2\,-2\,-1\,0:-2\,0$ & $2\,2\,1\,0.-1.2.2:2\,0$ \\ \hline
$K12n640$ & $3\,1:-2\,-1\,0:-2\,0$ & $3\,1\,0.-1.3.2:2\,0$ \\ \hline
$K12n644$ & $4:-2\,-1\,0:-2\,0$ & $2\,0.-1.2\,1.2:4\,0$ \\ \hline
$K12n647$ & $-2\,-1.2.-2.3\,0$ & $6^*2.2\,1.3.2\,0.2\,0.-1$ \\
\hline $K12n650$ & $8^*-2\,-1\,0.2\,0.-2\,0$ &
$102^*2\,1\,0.-1.2.1.1.1.1.1.1.1$ \\ \hline $K12n655$ &
$3.-2.2.-2\,-1$ & $111^*1.-1.1.1.1.1.2\,0.3\,0.2\,0.2\,0.1$ \\
\hline $K12n683$ & $2.-2.3.2.2$ & $2.4\,-1\,1\,1.3.2.2$ \\ \hline
$K12n694$ & $2:(-2,-3)\,0:-2\,0$ & $(2,3)\,0.-1.2.2:2\,0$ \\ \hline
$K12n702$ & $102^*.-2.2$ & $102^*.4\,-1\,1\,1.2$ \\ \hline $K12n703$
& $102^*.-2:2$ & $102^*.4\,-1\,1\,1:2$ \\ \hline $K12n704$ &
$9^*3\,0:.-2\,0$ & $9^*3\,0:.4\,-1\,1\,1\,0$ \\ \hline

\end{tabular}

\begin{tabular}{|c|c|c|} \hline
$K12n705$ & $101^*2:::.-2\,0$ & $101^*2:::.4\,-1\,1\,1\,0$ \\ \hline
$K12n707$ & $2.2.2.-2.3$ & $2.2.2.4\,-1\,1\,1.3$ \\ \hline $K12n708$
& $-2.2.3.2.2$ & $4\,-1\,1\,1.2.3.2.2$ \\ \hline $K12n712$ &
$9^*3:.-2\,0$ & $9^*3:.4\,-1\,1\,1\,0$ \\ \hline $K12n713$ &
$8^*3\,0:2:-2\,0$ & $8^*3\,0:2:4\,-1\,1\,1\,0$ \\ \hline $K12n715$ &
$8^*3\,0.2:-2\,0$ & $8^*3\,0.2:4\,-1\,1\,1\,0$ \\ \hline $K12n717$ &
$8^*-2.-2\,0.-2\,0:2$ & $9^*2\,0.2\,0.1.1.1.2\,0.1.-1.2\,0$ \\
\hline $K12n718$ & $8^*2.3.-2\,0$ & $8^*2.3.4\,-1\,1\,1\,0$ \\
\hline $K12n719$ & $2.-2.-2.-2.3\,0$ & $9^*2.2.1.1.1.3\,0.1.-1.1$ \\
\hline $K12n720$ & $8^*2\,0.-2\,0.2:2\,0$ &
$8^*2\,0.4\,-1\,1\,1\,0.2:2\,0$
\\ \hline $K12n728$ & $102^*2\,0::-2$ & $102^*2\,0::4\,-1\,1\,1$ \\
\hline $K12n732$ & $101^*2::-2\,0$ & $101^*2::4\,-1\,1\,1\,0$ \\
\hline $K12n734$ & $8^*2\,0.3\,0:.-2\,0$ &
$8^*2\,0.3\,0:.4\,-1\,1\,1\,0$ \\ \hline $K12n735$ &
$8^*2\,0.3:.-2\,0$ & $8^*2\,0.3:.4\,-1\,1\,1\,0$ \\ \hline $K12n737$
& $8^*-3\,-1\,0::2\,0$ & $8^*4\,-1::2\,0$ \\ \hline $K12n739$ &
$8^*2.-2.-2\,0.2\,0$ & $101^*2.-1.2\,0.2\,0.1.1.1.1.1.1$ \\ \hline
$K12n741$ & $8^*2:.3\,0:.-2\,0$ & $8^*2:.3\,0:.4\,-1\,1\,1\,0$ \\
\hline $K12n743$ & $9^*:3\,0.-2\,0$ & $9^*:3\,0.4\,-1\,1\,1\,0$ \\
\hline $K12n744$ & $2.3.-2.2\,0.2\,0$ & $2.3.4\,-1\,1\,1.2\,0.2\,0$
\\ \hline $K12n745$ & $8^*2.3\,0.-2\,0$ & $8^*2.3\,0.4\,-1\,1\,1\,0$
\\ \hline $K12n747$ & $2.-2.-2.-2.3$ & $9^*2.2.2.1.1.2.1.-1.1$ \\
\hline $K12n748$ & $2.2.2.2.-3$ & $2.2.2.2.2\,1\,-1\,0$ \\ \hline
$K12n759$ & $8^*3\,0.-2\,0:.2$ & $8^*3\,0.4\,-1\,1\,1\,0:.2$ \\
\hline $K12n764$ & $9^*2.-2.-2\,0$ & $8^*2.2.1.2.1.2\,0.2\,0.-1$ \\
\hline $K12n769$ & $8^*.2:2\,0.2\,0:-2\,0$ &
$8^*.2:2\,0.2\,0:4\,-1\,1\,1\,0$
\\ \hline $K12n783$ & $8^*-3\,0.2\,0.2\,0$ & $8^*2\,1\,-1.2\,0.2\,0$
\\ \hline $K12n787$ & $3.2.-2.-2\,0.-2$ &
$101^*2.1.1.2\,0.1.2\,0.1.2\,1\,0.1.-1$ \\ \hline $K12n788$ &
$-2.-2\,0.-2.3.2\,0$ & $111^*2.2\,0.2.2.-1.1.1.1.2.1.2\,0$ \\ \hline
$K12n792$ & $2.-4.-2\,0.-2$ & $8^*2.2\,0.1.2\,0.3.1.2.-1$ \\ \hline
$K12n794$ & $2.-3\,-1.-2\,0.-2$ & $8^*2.2\,0.1.2\,0.3\,0.1.2.-1$ \\
\hline $K12n798$ & $-2.-2.-2.2.2\,0.2\,0$ &
$9^*2\,0.2\,0.1.1.2\,0.2\,0.1.-1.1$ \\ \hline $K12n806$ &
$8^*-2\,0:-2\,0:-2\,0:-2\,0$ & $9^*2.-1.2\,0.2\,0.1.1.1.1.2$ \\
\hline $K12n808$ & $8^*-2\,0:-2\,0:2\,0:2\,0$ &
$8^*2.2\,0.1.2\,0.1.2\,0.2\,0.-1$ \\ \hline $K12n812$ &
$9^*.-2:-2\,0.2$ & $9^*2.2.1.-1.2.1.2\,0.1.1$ \\ \hline $K12n814$ &
$8^*2\,0:-3\,0:2\,0$ & $8^*2\,0:2\,1\,-1:2\,0$ \\ \hline $K12n815$ &
$9^*.-2:2\,0.-2$ & $8^*2\,0.2\,0.-1.2\,0.2\,0.1.2\,0.1$ \\ \hline
$K12n820$ & $3\,0.2.2\,0.-3\,0$ & $3\,0.2.2\,0.2\,1\,-1$ \\ \hline
$K12n821$ & $2:5\,0:-2\,0$ & $2:5\,0:4\,-1\,1\,1\,0$ \\ \hline
$K12n822$ & $-3\,0.2.2\,0.3\,0$ & $2\,1\,-1.2.2\,0.3\,0$ \\ \hline
$K12n823$ & $-3\,-2:2\,0:2\,0$ & $3\,1\,1\,-1\,0:2\,0:2\,0$ \\
\hline $K12n824$ & $2:3\,2\,0:-2\,0$ & $2:3\,2\,0:4\,-1\,1\,1\,0$ \\
\hline $K12n825$ & $-2.2.3.3\,0$ & $4\,-1\,1\,1.2.3.3\,0$ \\ \hline
$K12n829$ & $2.-2.-3.-2.2$ & $9^*2.3.1.1.1.2.1.-1.1$ \\ \hline
$K12n832$ & $2\,0.2.-3.-2\,0.-2$ & $9^*2.3\,0.1.1.1.2.1.-1.1$ \\
\hline $K12n837$ & $1212^*-1.-1.-1.-1.-1.-1$ &
$1310^*1.-1.1.1.1.1.1.1.1.1.1.1.1$ \\ \hline $K12n839$ &
$111^*:.-2\,0$ & $111^*:.4\,-1\,1\,1\,0$ \\ \hline $K12n840$ &
$101^*-2\,0::.2\,0$ & $101^*4\,-1\,1\,1\,0::.2\,0$ \\ \hline

\end{tabular}

\begin{tabular}{|c|c|c|} \hline
$K12n841$ & $8^*2.-2\,0.2\,0.2\,0$ & $8^*2.4\,-1\,1\,1\,0.2\,0.2\,0$
\\ \hline $K12n843$ & $9^*-2\,0:.-2\,0:.-2\,0$ &
$102^*2.2\,0.1.1.1.2\,0.-1.2\,0.1.1$ \\ \hline $K12n844$ &
$8^*2.2\,0.2:.-2\,0$ & $8^*2.2\,0.2:.4\,-1\,1\,1\,0$ \\ \hline
$K12n845$ & $9^*-2\,0::3$ & $9^*4\,-1\,1\,1\,0::3$ \\ \hline
$K12n846$ & $8^*4:-2\,0$ & $8^*4:4\,-1\,1\,1\,0$ \\ \hline $K12n848$
& $8^*3\,0.2\,0.-2\,0$ & $8^*3\,0.2\,0.4\,-1\,1\,1\,0$ \\ \hline
$K12n850$ & $4\,1:-2\,0:-2\,0$ & $2\,0.-1.2.2:4\,1\,0$ \\ \hline
$K12n851$ & $3\,1\,1:-2\,0:-2\,0$ & $3\,1\,1\,0.-1.2.2:2\,0$ \\
\hline $K12n858$ & $8^*-2\,0.2\,0.2\,0:2$ &
$8^*4\,-1\,1\,1\,0.2\,0.2\,0:2$ \\ \hline $K12n867$ &
$8^*-3\,0.2\,0:2$ & $8^*2\,1\,-1.2\,0:2$ \\ \hline $K12n868$ &
$9^*.-2:.-2:.2\,0$ & $111^*1.2\,0.1.2\,0.1.-1.1.1.1.1.2$ \\ \hline
$K12n870$ & $9^*2\,0.-2:.-2$ & $8^*2.2\,0.2\,0.-1.2\,0.2\,0.1.1$ \\
\hline $K12n881$ & $2.-2.-2\,0.-2.2.2\,0$ &
$9^*2.2\,0.1.1.2\,0.2.1.-1.1$ \\ \hline $K12n882$ &
$2\,0.3.2\,0.-3\,0$ & $2\,0.3.2\,0.2\,1\,-1$ \\ \hline

\end{tabular}

\bigskip

In the next table almost alternating knots from the first table are
given by their minimal DT-codes and DT-codes of their almost
alternating representations:

\begin{landscape}

\footnotesize

\begin{tabular}{|c|c|c|} \hline
$K12n176$ & $\{\{12\},\{4,8,14,2,-18,16,6,20,22,-24,12,-10\}\}$ &
$\{\{17\},\{4,8,14,2,24,32,6,30,26,28,-16,12,34,18,20,22,10\}\}$ \\
\hline $K12n204$ &
$\{\{12\},\{4,8,14,2,-18,-22,6,20,-10,24,-12,-16\}\}$ &
$\{\{17\},\{4,8,30,2,22,32,-18,24,26,28,34,14,16,12,6,20,10\}\}$ \\
\hline $K12n257$ &
$\{\{12\},\{4,8,16,2,-20,18,-22,6,12,-24,-14,-10\}\}$ &
$\{\{13\},\{4,8,18,2,22,-26,20,24,6,12,16,10,14\}\}$ \\
\hline $K12n258$ &
$\{\{12\},\{4,8,16,2,-20,18,-22,6,12,24,-14,-10\}\}$ &
$\{\{13\},\{4,14,26,22,18,6,20,2,10,24,12,8,-16\}\}$ \\
\hline $K12n281$ &
$\{\{12\},\{4,10,12,-14,16,2,-22,20,24,8,-6,18\}\}$ &
$\{\{17\},\{4,10,12,24,28,2,26,30,34,32,-16,8,14,6,22,20,18\}\}$ \\
\hline $K12n282$ &
$\{\{12\},\{4,10,12,14,-16,2,22,-20,24,-8,6,18\}\}$ &
$\{\{13\},\{4,10,26,16,20,2,18,22,6,24,8,14,-12\}\}$ \\
\hline $K12n286$ &
$\{\{12\},\{4,10,12,-14,20,2,-18,22,-6,24,16,8\}\}$ &
$\{\{17\},\{4,10,12,16,26,2,34,8,32,28,30,-18,14,6,20,22,24\}\}$ \\
\hline $K12n299$ &
$\{\{12\},\{4,10,12,-16,-20,2,18,-22,-6,24,-8,-14\}\}$ &
$\{\{15\},\{4,10,30,18,24,2,28,22,26,6,-12,16,8,14,20\}\}$ \\
\hline $K12n316$ &
$\{\{12\},\{4,10,-14,12,2,18,-20,22,24,8,-6,16\}\}$ &
$\{\{17\},\{4,14,24,22,20,30,2,28,6,10,8,-18,32,34,12,16,26\}\}$ \\
\hline $K12n317$ &
$\{\{12\},\{4,10,-14,12,2,18,20,-22,-24,8,-6,-16\}\}$ &
$\{\{15\},\{4,12,18,-14,26,2,30,20,6,16,8,28,10,22,24\}\}$ \\
\hline $K12n319$ &
$\{\{12\},\{4,10,14,-12,2,-18,20,22,24,-8,6,16\}\}$ &
$\{\{14\},\{4,12,24,-28,16,2,20,22,26,8,14,10,6,18\}\}$ \\
\hline $K12n321$ &
$\{\{12\},\{4,10,-14,12,2,18,-20,24,22,8,-6,16\}\}$ &
$\{\{17\},\{4,16,20,24,22,-6,30,2,28,12,10,8,34,32,14,18,26\}\}$ \\
\hline $K12n322$ &
$\{\{12\},\{4,10,-14,12,2,18,20,-24,-22,8,-6,-16\}\}$ &
$\{\{15\},\{4,12,18,-14,26,2,30,20,6,16,8,28,10,24,22\}\}$ \\
\hline $K12n324$ &
$\{\{12\},\{4,10,14,-12,2,-18,20,24,22,-8,6,16\}\}$ &
$\{\{15\},\{4,10,16,22,2,28,-24,6,30,26,14,8,20,12,18\}\}$ \\
\hline $K12n327$ &
$\{\{12\},\{4,10,-14,12,2,20,-18,22,-6,8,24,16\}\}$ &
$\{\{17\},\{4,10,26,22,2,14,24,34,32,30,6,12,8,16,20,18,-28\}\}$ \\
\hline $K12n349$ &
$\{\{12\},\{4,10,-14,-16,18,2,-20,-6,24,-22,-12,8\}\}$ &
$\{\{13\},\{4,10,-24,16,20,2,22,6,8,26,14,12,18\}\}$ \\
\hline $K12n353$ &
$\{\{12\},\{4,10,-14,16,18,2,22,-24,20,8,12,-6\}\}$ &
$\{\{13\},\{4,10,20,22,16,2,24,-12,26,6,18,14,8\}\}$ \\
\hline $K12n357$ &
$\{\{12\},\{4,10,-14,-16,20,2,18,-24,-22,12,8,-6\}\}$ &
$\{\{14\},\{4,16,12,-26,18,24,28,20,2,8,14,6,10,22\}\}$ \\
\hline $K12n359$ &
$\{\{12\},\{4,10,-14,16,20,2,18,-24,22,12,8,-6\}\}$ &
$\{\{15\},\{4,10,22,26,24,2,6,20,30,16,-28,12,8,14,18\}\}$ \\
\hline $K12n360$ &
$\{\{12\},\{4,10,14,-16,-20,2,18,24,-22,12,-8,6\}\}$ &
$\{\{13\},\{6,10,18,14,24,20,-2,26,22,12,4,16,8\}\}$ \\
\hline $K12n361$ &
$\{\{12\},\{4,10,14,-16,20,2,18,24,-22,12,8,-6\}\}$ &
$\{\{17\},\{4,14,18,30,28,26,24,2,22,34,12,16,20,-32,8,6,10\}\}$ \\
\hline $K12n376$ &
$\{\{12\},\{4,10,-14,18,2,20,-24,-6,-22,8,12,-16\}\}$ &
$\{\{14\},\{6,12,16,28,24,18,4,22,2,26,14,8,20,-10\}\}$ \\
\hline $K12n381$ &
$\{\{12\},\{4,10,14,-18,2,-22,8,-24,-20,-6,-12,-16\}\}$ &
$\{\{14\},\{4,10,24,12,18,2,20,26,-22,8,14,28,6,16\}\}$ \\
\hline $K12n382$ &
$\{\{12\},\{4,10,14,-18,2,22,8,-24,-20,-6,12,-16\}\}$ &
$\{\{13\},\{4,12,22,14,20,-16,2,8,26,10,24,6,18\}\}$ \\
\hline $K12n394$ &
$\{\{12\},\{4,10,-14,-18,-16,2,20,-24,-22,-6,12,-8\}\}$ &
$\{\{17\},\{4,14,26,32,30,18,20,2,10,28,34,8,6,16,-12,24,22\}\}$ \\
\hline $K12n395$ &
$\{\{12\},\{4,10,-14,-18,16,2,-20,24,22,-6,-12,8\}\}$ &
$\{\{13\},\{4,10,-12,22,16,2,20,18,8,26,24,6,14\}\}$ \\
\hline $K12n398$ &
$\{\{12\},\{4,10,-14,-18,-16,2,22,-24,-20,-8,12,-6\}\}$ &
$\{\{13\},\{4,10,24,20,-14,2,22,26,8,6,18,12,16\}\}$ \\
\hline $K12n399$ &
$\{\{12\},\{4,10,-14,18,16,2,22,-24,20,8,12,-6\}\}$ &
$\{\{13\},\{4,10,20,18,16,2,24,-12,26,22,6,14,8\}\}$ \\
\hline $K12n403$ &
$\{\{12\},\{4,10,-14,-18,22,2,-20,24,-6,-8,-12,16\}\}$ &
$\{\{17\},\{4,16,34,-12,30,20,28,24,2,32,10,14,18,6,8,22,26\}\}$ \\
\hline $K12n417$ &
$\{\{12\},\{4,10,-14,-20,2,18,-16,-6,24,-22,-8,12\}\}$ &
$\{\{13\},\{4,12,-24,16,22,2,20,6,14,26,8,10,18\}\}$ \\
\hline $K12n425$ &
$\{\{12\},\{4,10,-14,-22,2,18,-8,20,12,24,-6,16\}\}$ &
$\{\{13\},\{4,10,20,12,2,18,-24,22,8,26,14,16,6\}\}$ \\
\hline $K12n428$ &
$\{\{12\},\{4,10,14,22,2,18,8,-20,12,-24,6,-16\}\}$ &
$\{\{15\},\{4,12,30,26,18,2,-6,20,24,10,28,16,8,22,14\}\}$ \\
\hline $K12n430$ &
$\{\{12\},\{4,10,-14,-22,-16,2,20,-24,-8,-6,12,-18\}\}$ &
$\{\{16\},\{6,14,12,30,28,24,22,20,8,10,-2,4,32,18,16,26\}\}$ \\
\hline $K12n434$ &
$\{\{12\},\{4,10,14,-22,-20,2,-18,24,-8,-12,-16,-6\}\}$ &
$\{\{13\},\{4,10,12,20,22,2,-26,8,6,24,16,14,18\}\}$ \\
\hline $K12n436$ &
$\{\{12\},\{4,10,-14,-24,-18,2,22,-6,-20,-8,12,-16\}\}$ &
$\{\{13\},\{4,10,16,20,-14,2,22,26,24,8,18,12,6\}\}$ \\
\hline $K12n440$ &
$\{\{12\},\{4,10,-16,12,2,18,20,22,-24,8,14,-6\}\}$ &
$\{\{17\},\{4,10,26,12,2,22,24,6,34,30,32,8,14,16,18,20,-28\}\}$ \\
\hline $K12n444$ &
$\{\{12\},\{4,10,-16,12,2,18,22,20,-24,8,14,-6\}\}$ &
$\{\{13\},\{6,12,20,-22,16,24,2,10,8,4,26,14,18\}\}$ \\
\hline
\end{tabular}

\begin{tabular}{|c|c|c|} \hline

$K12n450$ & $\{\{12\},\{4,10,-16,-14,2,18,-22,-20,24,-8,-6,-12\}\}$
&
$\{\{17\},\{4,10,24,26,2,28,8,6,-22,30,32,34,14,12,18,20,16\}\}$ \\
\hline $K12n453$ &
$\{\{12\},\{4,10,16,-14,2,18,22,-20,24,8,-6,12\}\}$ &
$\{\{17\},\{4,14,20,28,26,24,2,22,32,12,34,18,-30,8,6,10,16\}\}$ \\
\hline $K12n459$ &
$\{\{12\},\{4,10,16,-14,2,24,22,-20,12,8,-6,18\}\}$ &
$\{\{17\},\{4,10,20,32,2,34,6,24,28,26,12,8,-30,18,16,14,22\}\}$ \\
\hline $K12n463$ &
$\{\{12\},\{4,10,16,-14,18,2,-22,20,12,24,-6,8\}\}$ &
$\{\{17\},\{4,10,22,20,24,2,6,30,28,26,12,34,8,-32,16,14,18\}\}$ \\
\hline $K12n465$ &
$\{\{12\},\{4,10,16,14,-20,2,-18,24,22,-8,-12,6\}\}$ &
$\{\{17\},\{-6,14,12,16,22,26,4,2,34,24,10,32,30,8,18,20,28\}\}$ \\
\hline $K12n481$ &
$\{\{12\},\{4,10,-16,18,14,2,20,24,-22,8,12,-6\}\}$ &
$\{\{15\},\{4,10,12,24,18,2,22,28,8,26,-30,6,16,20,14\}\}$ \\
\hline $K12n485$ &
$\{\{12\},\{4,10,-16,20,2,18,8,22,-24,12,14,-6\}\}$ &
$\{\{14\},\{6,22,14,16,-24,20,4,28,26,10,2,12,18,8\}\}$ \\
\hline $K12n486$ &
$\{\{12\},\{4,10,-16,-20,2,18,-22,-6,24,-8,-12,-14\}\}$ &
$\{\{14\},\{4,12,20,-28,16,2,22,24,26,6,10,14,8,18\}\}$ \\
\hline $K12n489$ &
$\{\{12\},\{4,10,16,20,2,18,-22,6,24,8,-12,-14\}\}$ &
$\{\{13\},\{4,10,16,20,2,18,22,6,26,8,24,14,-12\}\}$ \\
\hline $K12n497$ &
$\{\{12\},\{4,10,-16,-20,2,24,18,-6,22,-8,-14,12\}\}$ &
$\{\{13\},\{4,12,-14,16,24,20,2,18,6,26,10,8,22\}\}$ \\
\hline $K12n508$ &
$\{\{12\},\{4,10,-18,12,2,16,20,8,-24,22,14,-6\}\}$ &
$\{\{15\},\{4,12,24,22,18,2,28,-8,20,10,16,26,6,30,14\}\}$ \\
\hline $K12n509$ &
$\{\{12\},\{4,10,-18,12,16,2,-22,8,-24,-6,-14,-20\}\}$ &
$\{\{14\},\{8,10,22,16,2,24,20,28,6,26,12,4,-14,18\}\}$ \\
\hline $K12n510$ &
$\{\{12\},\{4,10,18,-12,-16,2,-22,-8,24,6,-14,-20\}\}$ &
$\{\{13\},\{4,10,20,-18,16,2,22,8,14,26,24,12,6\}\}$ \\
\hline $K12n519$ &
$\{\{12\},\{4,10,-18,-14,-16,2,20,-8,-24,22,12,-6\}\}$ &
$\{\{13\},\{4,10,24,16,-20,2,22,6,8,26,14,12,18\}\}$ \\
\hline $K12n520$ &
$\{\{12\},\{4,10,-18,-14,16,2,-22,20,12,-24,8,-6\}\}$ &
$\{\{13\},\{4,12,-24,22,6,16,2,20,26,10,14,8,18\}\}$ \\
\hline $K12n521$ &
$\{\{12\},\{4,10,18,-14,-16,2,-22,-20,-12,24,-8,-6\}\}$ &
$\{\{17\},\{4,10,26,24,16,2,6,8,34,32,30,14,12,18,22,20,-28\}\}$ \\
\hline $K12n524$ &
$\{\{12\},\{4,10,-18,16,14,2,20,8,-24,22,12,-6\}\}$ &
$\{\{14\},\{6,20,14,28,22,18,26,4,24,10,2,-12,16,8\}\}$ \\
\hline $K12n527$ &
$\{\{12\},\{4,10,18,-22,-16,2,-8,-20,-12,24,-14,-6\}\}$ &
$\{\{17\},\{4,10,26,12,16,2,24,8,32,30,34,14,6,22,20,18,-28\}\}$ \\
\hline $K12n535$ &
$\{\{12\},\{4,10,-20,16,14,2,-18,24,22,-12,-6,8\}\}$ &
$\{\{13\},\{8,14,22,16,-24,20,6,2,26,12,10,4,18\}\}$ \\
\hline $K12n536$ &
$\{\{12\},\{4,10,20,-16,-14,2,18,-24,-22,12,6,-8\}\}$ &
$\{\{13\},\{4,12,22,18,24,20,2,10,8,-26,14,6,16\}\}$ \\
\hline $K12n542$ &
$\{\{12\},\{4,10,20,-24,-14,2,18,-8,-22,12,6,-16\}\}$ &
$\{\{13\},\{4,12,22,18,16,20,2,10,24,-26,14,6,8\}\}$ \\
\hline $K12n544$ &
$\{\{12\},\{4,10,-22,-14,2,18,-8,-20,24,12,-6,16\}\}$ &
$\{\{13\},\{4,10,22,14,2,18,-24,6,26,12,16,8,20\}\}$ \\
\hline $K12n551$ &
$\{\{12\},\{4,12,14,18,-16,20,2,24,-22,6,10,-8\}\}$ &
$\{\{13\},\{-4,12,16,22,24,2,18,6,26,10,8,20,14\}\}$ \\
\hline $K12n552$ &
$\{\{12\},\{4,12,14,-18,24,-20,2,8,22,-6,-10,16\}\}$ &
$\{\{13\},\{4,14,-22,18,16,20,26,2,24,6,12,8,10\}\}$ \\
\hline $K12n557$ &
$\{\{12\},\{4,12,-16,14,20,2,8,-24,22,10,18,-6\}\}$ &
$\{\{13\},\{4,14,20,26,-12,24,2,8,22,6,18,10,16\}\}$ \\
\hline $K12n559$ &
$\{\{12\},\{4,12,-16,14,20,2,18,22,-24,8,10,-6\}\}$ &
$\{\{17\},\{6,10,22,30,32,24,28,26,-12,4,34,18,16,14,2,8,20\}\}$ \\
\hline $K12n561$ &
$\{\{12\},\{4,12,-16,14,20,2,18,22,-24,10,8,-6\}\}$ &
$\{\{17\},\{6,10,22,32,30,24,28,26,-12,4,34,18,16,14,2,8,20\}\}$ \\
\hline $K12n562$ &
$\{\{12\},\{4,12,16,-14,-20,2,-18,22,24,-10,-8,6\}\}$ &
$\{\{17\},\{4,12,26,16,14,2,24,8,6,34,30,32,10,18,20,22,-28\}\}$ \\
\hline $K12n563$ &
$\{\{12\},\{4,12,16,-14,-20,2,18,-22,24,-10,-8,6\}\}$ &
$\{\{17\},\{4,12,34,32,22,2,26,28,24,6,10,8,-30,14,16,18,20\}\}$ \\
\hline $K12n567$ &
$\{\{12\},\{4,12,-16,18,14,2,20,22,-24,8,10,-6\}\}$ &
$\{\{17\},\{4,18,22,24,26,-6,34,20,2,30,32,12,8,10,14,16,28\}\}$ \\
\hline $K12n568$ &
$\{\{12\},\{4,12,-16,18,20,2,8,22,-24,14,10,-6\}\}$ &
$\{\{17\},\{6,14,26,12,24,34,2,8,-22,32,30,28,4,10,16,20,18\}\}$ \\
\hline $K12n569$ &
$\{\{12\},\{4,12,-16,18,20,2,10,22,-24,8,14,-6\}\}$ &
$\{\{14\},\{4,14,28,22,-18,26,2,8,24,12,6,16,10,20\}\}$ \\
\hline $K12n579$ &
$\{\{12\},\{4,12,16,-20,14,-22,2,8,24,-6,-10,18\}\}$ &
$\{\{13\},\{4,14,18,22,16,24,20,2,8,26,-10,6,12\}\}$ \\
\hline $K12n580$ &
$\{\{12\},\{4,12,-16,20,18,2,10,22,-24,8,14,-6\}\}$ &
$\{\{14\},\{4,14,20,26,16,22,2,24,-6,28,12,10,8,18\}\}$ \\
\hline $K12n584$ &
$\{\{12\},\{4,12,-16,-22,14,2,18,20,-24,10,-8,-6\}\}$ &
$\{\{17\},\{4,12,26,24,14,2,22,8,34,30,32,10,6,16,18,20,-28\}\}$ \\
\hline $K12n589$ &
$\{\{12\},\{4,12,18,-14,16,2,20,-22,10,24,-8,6\}\}$ &
$\{\{17\},\{4,18,26,28,30,32,-8,20,2,24,6,34,16,22,14,10,12\}\}$ \\
\hline $K12n594$ &
$\{\{12\},\{4,12,-18,14,-20,2,8,-24,-22,-10,-6,-16\}\}$ &
$\{\{13\},\{4,14,18,24,16,22,2,8,-12,26,6,10,20\}\}$ \\
\hline $K12n612$ &
$\{\{12\},\{4,12,20,-16,18,14,2,24,-22,10,6,-8\}\}$ &
$\{\{13\},\{-4,14,20,12,24,8,22,2,6,26,16,10,18\}\}$ \\
\hline $K12n617$ &
$\{\{12\},\{4,14,-10,-18,-16,-24,20,2,-22,12,-8,-6\}\}$ &
$\{\{16\},\{6,14,12,24,16,20,22,28,30,32,4,2,18,10,-8,26\}\}$ \\
\hline $K12n618$ &
$\{\{12\},\{4,14,10,-18,-16,24,-20,2,-22,-12,-8,-6\}\}$ &
$\{\{17\},\{4,18,30,28,26,34,32,22,12,2,14,6,10,8,-24,20,16\}\}$ \\
\hline $K12n619$ &
$\{\{12\},\{4,14,10,-18,16,24,20,2,-22,12,-8,6\}\}$ &
$\{\{17\},\{4,18,22,32,30,28,20,6,24,2,14,34,8,12,10,-26,16\}\}$ \\
\hline
\end{tabular}

\begin{tabular}{|c|c|c|} \hline
$K12n624$ & $\{\{12\},\{4,14,-10,-18,-24,20,2,-8,-22,12,-6,-16\}\}$
&
$\{\{13\},\{-6,10,20,24,2,18,22,26,12,4,16,8,14\}\}$ \\
\hline $K12n629$ &
$\{\{12\},\{4,14,-10,-22,-16,-18,20,2,-6,24,-8,12\}\}$ &
$\{\{15\},\{6,10,24,-22,18,4,20,28,8,26,30,2,12,16,14\}\}$ \\
\hline $K12n630$ &
$\{\{12\},\{4,14,-10,22,-16,-18,-20,2,-6,-24,8,-12\}\}$ &
$\{\{13\},\{-6,10,16,20,2,22,26,24,4,12,8,18,14\}\}$ \\
\hline $K12n631$ &
$\{\{12\},\{4,14,-10,22,-16,18,20,2,-6,24,8,12\}\}$ &
$\{\{17\},\{4,18,10,32,20,34,30,26,28,2,6,12,14,16,-24,8,22\}\}$ \\
\hline $K12n638$ &
$\{\{12\},\{4,14,16,-18,24,-20,2,12,-22,-10,-6,8\}\}$ &
$\{\{13\},\{4,16,26,20,24,18,-22,2,8,12,6,10,14\}\}$ \\
\hline $K12n640$ &
$\{\{12\},\{4,14,-18,16,-20,-22,2,8,-24,-12,-10,-6\}\}$ &
$\{\{13\},\{8,-12,18,22,2,20,26,24,4,6,10,14,16\}\}$ \\
\hline $K12n644$ &
$\{\{12\},\{4,14,-18,16,-22,-20,2,8,-24,-12,-10,-6\}\}$ &
$\{\{13\},\{4,16,20,26,18,24,22,2,8,-14,6,12,10\}\}$ \\
\hline $K12n647$ &
$\{\{12\},\{4,14,-18,16,24,-20,2,10,-22,-12,-6,8\}\}$ &
$\{\{13\},\{4,16,20,24,26,18,22,2,10,-14,6,12,8\}\}$ \\
\hline $K12n650$ &
$\{\{12\},\{4,14,-18,-20,24,-6,2,-22,-12,-10,-16,8\}\}$ &
$\{\{13\},\{4,16,-20,12,26,24,22,2,8,10,14,6,18\}\}$ \\
\hline $K12n655$ &
$\{\{12\},\{4,14,-18,-22,16,-20,2,10,-24,-12,-8,-6\}\}$ &
$\{\{16\},\{6,10,12,22,18,16,26,4,-2,28,8,32,30,14,20,24\}\}$ \\
\hline $K12n683$ &
$\{\{12\},\{6,8,-12,16,-18,-20,22,24,2,-10,-4,14\}\}$ &
$\{\{17\},\{-6,16,18,14,28,22,34,2,4,10,30,12,32,8,20,24,26\}\}$ \\
\hline $K12n694$ &
$\{\{12\},\{6,8,16,2,-20,18,-22,4,12,-24,-14,-10\}\}$ &
$\{\{13\},\{6,8,18,2,22,-26,20,24,4,12,16,10,14\}\}$ \\
\hline $K12n702$ &
$\{\{12\},\{6,-10,-12,14,-18,-20,24,8,-22,-2,-4,-16\}\}$ &
$\{\{17\},\{6,10,12,14,24,26,34,8,-22,30,28,32,2,4,20,18,16\}\}$ \\
\hline $K12n703$ &
$\{\{12\},\{6,10,12,-14,18,20,-24,-8,-22,2,4,-16\}\}$ &
$\{\{17\},\{-6,14,12,16,28,30,4,2,34,26,10,8,18,32,22,20,24\}\}$ \\
\hline $K12n704$ &
$\{\{12\},\{6,-10,-12,14,-20,-18,24,-22,-2,-4,-16,-8\}\}$ &
$\{\{17\},\{-6,16,14,12,22,24,34,4,2,26,32,30,20,18,8,10,28\}\}$ \\
\hline $K12n705$ &
$\{\{12\},\{6,-10,12,14,20,-18,24,22,2,-4,16,8\}\}$ &
$\{\{17\},\{-6,18,20,16,32,24,22,34,2,4,28,30,8,14,12,10,26\}\}$ \\
\hline $K12n707$ &
$\{\{12\},\{6,-10,-12,16,18,-20,22,24,8,-2,-4,14\}\}$ &
$\{\{17\},\{-6,18,20,16,28,22,30,34,2,4,10,14,32,8,12,24,26\}\}$ \\
\hline $K12n708$ &
$\{\{12\},\{6,-10,12,-16,-18,20,22,24,-8,-2,4,14\}\}$ &
$\{\{17\},\{6,10,26,20,2,22,24,34,30,32,8,12,4,14,16,18,-28\}\}$ \\
\hline $K12n712$ &
$\{\{12\},\{6,-10,-12,16,-20,-18,22,24,-2,-4,-8,14\}\}$ &
$\{\{17\},\{-6,16,14,12,22,24,34,4,2,26,32,30,18,20,8,10,28\}\}$ \\
\hline $K12n713$ &
$\{\{12\},\{6,-10,12,16,20,-18,22,24,2,-4,8,14\}\}$ &
$\{\{17\},\{-6,16,18,20,32,24,30,2,4,34,28,10,12,8,14,22,26\}\}$ \\
\hline $K12n715$ &
$\{\{12\},\{6,-10,12,16,20,-18,24,22,2,-4,8,14\}\}$ &
$\{\{17\},\{-6,14,16,18,32,28,2,4,34,26,30,10,8,12,22,20,24\}\}$ \\
\hline $K12n717$ &
$\{\{12\},\{6,-10,-12,16,22,-18,20,24,-2,-4,14,8\}\}$ &
$\{\{13\},\{-6,12,20,16,26,22,2,24,8,10,4,14,18\}\}$ \\
\hline $K12n718$ &
$\{\{12\},\{6,-10,12,16,22,-18,20,24,2,-4,14,8\}\}$ &
$\{\{17\},\{-6,14,16,18,32,30,2,4,34,26,28,10,8,20,12,22,24\}\}$ \\
\hline $K12n719$ &
$\{\{12\},\{6,10,12,16,-22,18,-20,-24,2,4,-14,-8\}\}$ &
$\{\{13\},\{6,12,18,26,-16,22,4,20,2,24,10,14,8\}\}$ \\
\hline $K12n720$ &
$\{\{12\},\{6,-10,-12,16,22,-20,18,8,24,-4,-2,14\}\}$ &
$\{\{17\},\{6,18,10,22,2,26,28,24,4,32,8,34,-30,12,14,16,20\}\}$ \\
\hline $K12n728$ &
$\{\{12\},\{6,-10,-12,18,-22,-16,-20,-4,-2,24,-14,-8\}\}$ &
$\{\{17\},\{6,10,16,26,24,4,8,2,32,30,34,12,14,22,20,18,-28\}\}$ \\
\hline $K12n732$ &
$\{\{12\},\{6,-10,12,22,20,-18,24,8,2,-4,16,14\}\}$ &
$\{\{17\},\{-6,18,20,16,32,30,22,34,2,4,28,10,8,14,12,24,26\}\}$ \\
\hline $K12n734$ &
$\{\{12\},\{6,-10,12,24,20,-18,8,22,2,-4,14,16\}\}$ &
$\{\{17\},\{6,10,26,14,4,22,24,34,32,30,8,12,2,16,20,18,-28\}\}$ \\
\hline $K12n735$ &
$\{\{12\},\{6,-10,12,24,20,-18,8,22,2,-4,16,14\}\}$ &
$\{\{17\},\{6,10,26,14,4,24,22,34,32,30,8,12,2,16,20,18,-28\}\}$ \\
\hline $K12n737$ &
$\{\{12\},\{6,-10,12,24,22,-18,20,8,2,-4,14,16\}\}$ &
$\{\{13\},\{6,10,20,12,4,18,22,24,8,2,-26,14,16\}\}$ \\
\hline $K12n739$ &
$\{\{12\},\{6,10,14,16,-20,4,-22,2,24,-8,-12,18\}\}$ &
$\{\{13\},\{6,12,16,18,20,-24,4,22,2,26,10,14,8\}\}$ \\
\hline $K12n741$ &
$\{\{12\},\{6,-10,-14,16,-22,-18,-20,24,-2,-12,-4,-8\}\}$ &
$\{\{17\},\{-6,18,20,16,26,22,24,34,2,4,30,12,32,8,14,10,28\}\}$ \\
\hline $K12n743$ &
$\{\{12\},\{6,-10,-14,16,-22,-20,-18,24,-2,-12,-4,-8\}\}$ &
$\{\{17\},\{6,14,18,28,32,30,26,24,22,12,2,16,4,20,-34,10,8\}\}$ \\
\hline $K12n744$ &
$\{\{12\},\{6,-10,14,16,-22,20,18,24,2,12,-4,-8\}\}$ &
$\{\{17\},\{-6,18,16,20,30,26,22,28,4,2,34,12,32,14,10,8,24\}\}$ \\
\hline $K12n745$ &
$\{\{12\},\{6,-10,14,16,22,-20,-18,24,2,-12,-4,8\}\}$ &
$\{\{17\},\{-6,14,16,18,32,30,2,4,34,28,26,10,8,20,12,22,24\}\}$ \\
\hline $K12n747$ &
$\{\{12\},\{6,-10,-14,18,-2,-20,-4,22,24,8,-12,16\}\}$ &
$\{\{13\},\{6,10,14,20,-16,4,24,2,22,26,8,12,18\}\}$ \\
\hline $K12n748$ &
$\{\{12\},\{6,10,-14,18,2,-20,-4,-22,-24,8,-12,-16\}\}$ &
$\{\{13\},\{4,10,18,22,2,-26,20,24,12,6,14,8,16\}\}$ \\
\hline $K12n759$ &
$\{\{12\},\{6,-10,14,24,22,-20,-18,8,2,-12,-4,16\}\}$ &
$\{\{17\},\{6,10,22,34,26,24,4,-20,28,30,32,12,2,16,18,14,8\}\}$ \\
\hline $K12n764$ &
$\{\{12\},\{6,-10,-16,14,22,-2,20,24,-4,8,12,18\}\}$ &
$\{\{13\},\{8,14,18,22,26,20,4,24,2,12,6,-10,16\}\}$ \\
\hline
\end{tabular}

\begin{tabular}{|c|c|c|} \hline
$K12n769$ &
$\{\{12\},\{6,-10,16,-18,-22,-4,-20,-24,2,-14,-8,-12\}\}$ &
$\{\{17\},\{-6,20,22,18,32,28,24,30,34,2,4,12,16,8,14,10,26\}\}$ \\
\hline $K12n783$ &
$\{\{12\},\{6,-10,16,-22,-20,-2,-8,4,24,-14,-12,18\}\}$ &
$\{\{13\},\{4,14,18,20,24,6,2,-26,22,12,10,16,8\}\}$ \\
\hline $K12n787$ &
$\{\{12\},\{6,-10,-18,14,-2,-20,8,22,24,-4,-12,16\}\}$ &
$\{\{15\},\{4,16,20,12,28,22,6,2,26,30,-14,8,18,10,24\}\}$ \\
\hline $K12n788$ &
$\{\{12\},\{6,-10,-18,14,-2,-20,8,24,22,-4,-12,16\}\}$ &
$\{\{17\},\{6,14,26,34,32,18,28,2,22,30,10,16,-8,4,12,20,24\}\}$ \\
\hline $K12n792$ &
$\{\{12\},\{6,-10,-18,16,-2,-20,-22,8,24,-4,-12,-14\}\}$ &
$\{\{14\},\{-6,12,18,26,22,16,2,24,4,28,8,10,14,20\}\}$ \\
\hline $K12n794$ &
$\{\{12\},\{6,-10,-18,16,-2,-20,-22,8,24,-4,-14,-12\}\}$ &
$\{\{14\},\{-6,12,18,26,22,16,2,24,4,28,10,8,14,20\}\}$ \\
\hline $K12n798$ &
$\{\{12\},\{6,-10,18,-16,-22,20,4,-8,24,12,-2,-14\}\}$ &
$\{\{13\},\{8,14,22,16,24,4,20,2,26,-6,12,10,18\}\}$ \\
\hline $K12n806$ &
$\{\{12\},\{6,-10,-20,12,-16,-2,18,-22,-8,24,-4,-14\}\}$ &
$\{\{13\},\{8,12,18,26,22,4,20,-2,24,10,14,6,16\}\}$ \\
\hline $K12n808$ &
$\{\{12\},\{6,-10,-20,12,16,-2,18,22,8,24,-4,14\}\}$ &
$\{\{13\},\{-6,10,14,24,2,16,20,4,22,26,12,8,18\}\}$ \\
\hline $K12n812$ &
$\{\{12\},\{6,-10,-20,12,16,-2,24,22,8,14,-4,18\}\}$ &
$\{\{13\},\{-6,10,14,24,2,22,20,4,12,26,16,8,18\}\}$ \\
\hline $K12n814$ &
$\{\{12\},\{6,10,-20,12,-16,2,-24,-22,-8,-14,-4,-18\}\}$ &
$\{\{13\},\{4,12,-14,22,26,18,2,8,24,10,16,6,20\}\}$ \\
\hline $K12n815$ &
$\{\{12\},\{6,10,-20,12,-16,2,24,22,-8,14,-4,18\}\}$ &
$\{\{13\},\{6,10,22,14,2,-16,20,8,24,12,26,4,18\}\}$ \\
\hline $K12n820$ &
$\{\{12\},\{6,-10,-20,16,-2,-18,22,24,8,-4,12,14\}\}$ &
$\{\{13\},\{4,12,20,14,18,2,-26,22,24,10,6,16,8\}\}$ \\
\hline $K12n821$ &
$\{\{12\},\{6,-10,20,-16,-2,-18,22,24,-8,4,12,14\}\}$ &
$\{\{17\},\{6,12,22,24,26,28,4,-20,30,32,34,2,8,10,16,18,14\}\}$ \\
\hline $K12n822$ &
$\{\{12\},\{6,10,-20,16,2,-18,-22,-24,8,-4,-12,-14\}\}$ &
$\{\{13\},\{4,12,-10,22,18,2,20,24,8,26,6,14,16\}\}$ \\
\hline $K12n823$ &
$\{\{12\},\{6,-10,-20,16,-2,-18,22,24,8,-4,14,12\}\}$ &
$\{\{13\},\{6,14,16,22,18,20,4,2,-26,24,12,8,10\}\}$ \\
\hline $K12n824$ &
$\{\{12\},\{6,-10,20,-16,-2,-18,22,24,-8,4,14,12\}\}$ &
$\{\{17\},\{6,12,22,24,28,26,4,-20,30,32,34,2,8,10,16,18,14\}\}$ \\
\hline $K12n825$ &
$\{\{12\},\{6,10,-20,16,2,-18,-22,-24,8,-4,-14,-12\}\}$ &
$\{\{17\},\{-6,18,20,16,24,22,30,34,2,4,28,10,32,8,12,14,26\}\}$ \\
\hline $K12n829$ &
$\{\{12\},\{6,-10,-22,16,-2,-18,-20,24,8,-12,-4,14\}\}$ &
$\{\{13\},\{6,10,22,16,-24,4,20,26,8,12,14,2,18\}\}$ \\
\hline $K12n832$ &
$\{\{12\},\{6,-10,-22,16,-2,-20,-18,24,8,-12,-4,14\}\}$ &
$\{\{13\},\{6,10,22,16,-24,4,20,26,8,14,12,2,18\}\}$ \\
\hline $K12n837$ &
$\{\{12\},\{6,-12,-10,16,14,-20,-18,24,22,-4,-2,8\}\}$ &
$\{\{13\},\{6,20,22,18,16,4,2,24,26,14,12,-10,8\}\}$ \\
\hline $K12n839$ &
$\{\{12\},\{6,-12,10,16,14,-20,18,24,22,-4,2,8\}\}$ &
$\{\{17\},\{6,12,10,26,14,24,22,32,30,34,8,4,2,20,18,16,-28\}\}$ \\
\hline $K12n840$ &
$\{\{12\},\{6,-12,10,16,24,-20,18,8,22,-4,2,14\}\}$ &
$\{\{17\},\{6,12,24,26,14,4,22,32,30,34,8,10,2,20,18,16,-28\}\}$ \\
\hline $K12n841$ &
$\{\{12\},\{6,-12,-10,18,20,-16,-22,-4,-2,24,-14,8\}\}$ &
$\{\{17\},\{-6,20,22,18,26,32,28,30,34,2,4,14,12,8,24,16,10\}\}$ \\
\hline $K12n843$ &
$\{\{12\},\{6,-12,-10,24,14,-20,-18,8,22,-4,-2,16\}\}$ &
$\{\{14\},\{6,12,16,28,20,4,24,-10,2,26,8,18,14,22\}\}$ \\
\hline $K12n844$ &
$\{\{12\},\{6,-12,10,24,14,-20,18,8,22,-4,2,16\}\}$ &
$\{\{17\},\{-6,18,16,14,32,22,26,34,4,2,12,30,20,8,24,10,28\}\}$ \\
\hline $K12n845$ &
$\{\{12\},\{6,-12,-14,16,-20,-22,-18,24,-2,-4,-10,-8\}\}$ &
$\{\{17\},\{-6,16,14,12,22,24,34,4,2,30,32,28,18,8,10,20,26\}\}$ \\
\hline $K12n846$ &
$\{\{12\},\{6,-12,-14,16,-22,-20,-18,24,-2,-4,-10,-8\}\}$ &
$\{\{17\},\{-6,16,18,20,32,30,28,2,4,34,26,10,8,14,12,22,24\}\}$ \\
\hline $K12n848$ &
$\{\{12\},\{6,-12,-14,24,-20,-22,-18,8,-2,-4,-10,16\}\}$ &
$\{\{17\},\{-6,16,14,18,32,28,30,4,2,34,26,8,10,20,12,24,22\}\}$ \\
\hline $K12n850$ &
$\{\{12\},\{6,-12,-16,14,20,22,-2,24,-4,8,10,18\}\}$ &
$\{\{13\},\{6,16,18,22,26,20,24,2,4,10,-14,8,12\}\}$ \\
\hline $K12n851$ &
$\{\{12\},\{6,-12,-16,14,20,22,-2,24,-4,10,8,18\}\}$ &
$\{\{13\},\{6,16,18,22,-26,20,24,4,2,10,14,8,12\}\}$ \\
\hline $K12n858$ &
$\{\{12\},\{6,-12,16,24,-20,-4,-22,-8,2,-14,-10,-18\}\}$ &
$\{\{17\},\{-6,20,22,18,26,32,28,24,34,2,4,30,12,8,14,16,10\}\}$ \\
\hline $K12n867$ &
$\{\{12\},\{6,-12,-20,16,22,-2,-18,24,10,-4,14,8\}\}$ &
$\{\{13\},\{4,16,-14,20,24,26,22,2,10,8,12,6,18\}\}$ \\
\hline $K12n868$ &
$\{\{12\},\{6,-12,22,20,14,-4,18,-24,10,8,2,-16\}\}$ &
$\{\{14\},\{6,10,22,20,16,4,24,8,26,14,2,28,18,-12\}\}$ \\
\hline $K12n870$ &
$\{\{12\},\{6,-14,10,-20,2,24,18,-4,22,-8,12,16\}\}$ &
$\{\{13\},\{-6,10,24,18,2,20,26,22,12,8,16,4,14\}\}$ \\
\hline $K12n881$ &
$\{\{12\},\{6,-20,-10,24,14,-4,-18,8,22,-12,-2,16\}\}$ &
$\{\{13\},\{6,12,22,26,16,4,20,24,8,14,2,-10,18\}\}$ \\
\hline $K12n882$ &
$\{\{12\},\{6,-20,-12,16,18,-4,22,24,8,-2,-10,14\}\}$ &
$\{\{13\},\{4,12,-10,16,24,2,20,22,8,26,14,6,18\}\}$ \\ \hline

\end{tabular}
\end{landscape}

\normalsize

\bigskip

For the remaining 35 non-alternating knots it is not known that they
are almost alternating or not, so we don't know that their Turaev
genus is 1 or greater than 1. These knots are given in the following
table:

\bigskip

\begin{tabular}{|c|c|c|} \hline
$K12n253$ & $2:(2,-2\,-1)\,0:-2\,0$ \\ \hline $K12n254$ &
$2:(-2,2\,1)\,0:-2\,0$ \\ \hline $K12n280$ & $2.-2\,0.-2.2\,1\,1\,0$
\\ \hline $K12n285$ & $8^*-2\,-1\,-1::-2\,0$ \\ \hline $K12n318$ &
$2\,0.-2\,-1.-2\,0.3\,0$ \\ \hline $K12n323$ & $3.-2\,-1\,0.-2.2\,0$
\\ \hline $K12n328$ & $2\,0.-2\,-1.-2\,0.2\,1\,0$ \\ \hline
$K12n356$ & $8^*-2\,-1.2\,0.-2$ \\ \hline $K12n358$ &
$2.-3\,0.-2\,-1.2\,0$ \\ \hline $K12n370$ &
$-2\,-1\,0.-2\,-1.-2\,0.2\,0$ \\ \hline $K12n371$ &
$2\,1\,0.-2\,-1.-2\,0.2\,0$ \\ \hline $K12n375$ &
$2.2.-2.2\,0.-2\,-1$ \\ \hline $K12n407$ & $8^*2:.-2\,0:.-2\,-1\,0$
\\ \hline $K12n426$ & $8^*-2\,-1::-2\,-1$ \\ \hline $K12n439$ &
$8^*-3::-2\,-1$ \\ \hline $K12n443$ & $8^*-2\,-1::-3\,0$ \\ \hline
$K12n451$ & $8^*2:.-2\,-1\,0:.-2\,0$ \\ \hline $K12n452$ &
$3.-2\,0.-2.2\,1\,0$ \\ \hline $K12n462$ & $8^*-2\,0.2:-2\,-1\,0$ \\
\hline $K12n487$ & $2\,1\,0.-2.-2\,0.3\,0$ \\ \hline $K12n488$ &
$-2\,-1\,0.-2.-2\,0.3\,0$ \\ \hline $K12n591$ &
$3\,0.-2\,-1.-2\,0.2\,0$ \\ \hline $K12n706$ & $101^*-2\,0::.-2\,0$
\\ \hline $K12n729$ & $102^*-2\,0::-2$ \\ \hline $K12n730$ &
$8^*-2.2.-2\,0:2\,0$ \\ \hline $K12n749$ & $-2.2.-2.2.3$ \\ \hline
$K12n750$ & $2.2.-2.2.-3$ \\ \hline $K12n801$ &
$8^*2\,0.-2.-1.-2\,0.2\,0$ \\ \hline $K12n807$ &
$8^*-2\,0:-2\,0:-2\,0:2\,0$ \\ \hline $K12n809$ &
$8^*-2\,0:2\,0:-2\,0:2\,0$ \\ \hline $K12n811$ & $9^*.-2:-2\,0.-2$ \\
\hline $K12n830$ & $-2.2.-3.2.2$ \\ \hline $K12n835$ &
$8^*2\,0.-2\,0.-1.-2\,0.2\,0$ \\ \hline $K12n838$ &
$-2.-2.-2\,0.2.2.2\,0$ \\ \hline $K12n873$ &
$8^*2\,0.-2\,0.-2\,0.2\,0$ \\ \hline
\end{tabular}

\bigskip


\bigskip
\bigskip

\footnotesize

\noindent THE MATHEMATICAL INSTITUTE, KNEZ MIHAILOVA  36, P.O.BOX
367, \\ 11001 BELGRADE, SERBIA

\medskip

\noindent {\it E-mail address:} $\mathrm{sjablan@gmail.com}$

\bigskip

\end{document}